\documentclass{article}
\usepackage{fullpage}
\usepackage{amssymb,amsmath,amsfonts, amscd}
\usepackage{array}
\usepackage{amsmath}
\usepackage{amssymb}
\usepackage{titlefoot}

\begin{document}

\newtheorem{tha}{Theorem}
\renewcommand{\thetha}{\Alph{tha}}

\newtheorem{example}{Example}
\newtheorem{remark}{Remark}
\newtheorem{observation}{Observation}
\newtheorem{definition}{Definition}
\newtheorem{theorem}{Theorem} [section]
\newtheorem{corollary}[theorem]{Corollary}
\newtheorem{lemma}[theorem]{Lemma}
\newtheorem{fact}[theorem]{Fact}
\newtheorem{claim}{Claim}

\setlength{\textwidth}{17cm} \setlength{\oddsidemargin}{-0.1 in}
\setlength{\evensidemargin}{-0.1 in} \setlength{\topmargin}{0.0
in}

\def \df {\noindent {\bf Definition. }}

\newcommand{\dist}{{\rm dist}}
\newcommand{\Dist}{{\rm Dist}}
\newcommand{\chib}{\chi_B}
\newcommand{\Forb}{{\rm Forb}}

\newcommand{\qedbox}{$\blacksquare$ \newline}
\newenvironment{proof}%
{%
\noindent{\it Proof.} } {
\hfill\qedbox }

\newcommand{\proofend}{\hfill\qedbox}

\def\qed{\hskip 1.3em\hfill\rule{6pt}{6pt} \vskip 20pt}

\linespread{1.0}
\input epsf
\def\epsfsize#1#2{0.5#1\relax}
\def\O{\text{O}}
\def\o{\text{o}}
\def\ex{\text{ex}}
\def\Z{\mathbb Z}

\def\fl#1{\lfloor #1 \rfloor}
\def\ce#1{\lceil #1 \rceil}

\def \cH{{\cal H }}
\def \cF{{\cal F}}
\def \cG{{\cal G}}
\def \cQ{{\cal Q}}
\def \cA{{\cal A}}
\def \cD{{\cal D}}

\def \f {{\cal F }}
\def \A {{\cal A}}
\def \D {{\cal D}}
\def \fn2 {{\lfloor n/2 \rfloor}}
\def \cn2 {{\lceil  n/2 \rceil}}

\newcommand{\bin}[2]{{#1\choose #2}}
\newcommand{\comp}{\overline}
\newcommand{\exval}{{\rm ex}}

\newcommand{\subjclassname}{%
  \textup{2000} Mathematics Subject Classification}
\title{Avoiding patterns in matrices via a small
number of changes}\runningtitle{Avoiding patterns in matrices}
\author{Maria Axenovich\footnote{\tt axenovic@math.iastate.edu}
\hspace{.5in} Ryan Martin\footnote{\tt rymartin@iastate.edu} \\
Department of Mathematics \\ Iowa State University \\ Ames, IA
50011 \\ Phone: (515) 294-1752; Fax: (515) 294-5454}
 \amssubj{15A99, 05C15, 05B20, 05C50, 05C80}
\date{}
\maketitle

\begin{abstract}
   Let ${\cal A}=\{A_1,\ldots,  A_r\}$ be a partition of a set
   $\{1,\ldots,m\}\times\{1,\ldots, n\}$ into $r$ nonempty subsets,
   and $A=(a_{ij})$ be an $m\times n$ matrix.  We say that $A$ has
   a pattern ${\cal A}$ provided that  $a_{ij}=a_{i'j'}$ if and only if
   $(i,j),(i',j')\in A_t$ for some $t\in\{1,\ldots,r\}$.   In
   this note we study the following function $f$ defined on
   the set of all $m\times n$ matrices $M$ with $s$ distinct
   entries: $f(M; {\cal A})$ is the smallest number of positions
   where the entries of $M$  need to be changed such that the
   resulting matrix does not have any submatrix with pattern
   ${\cal A}$.  We give an asymptotically tight value for
   $$f(m,n; s, {\cal A}) = \max
   \{f(M; {\cal A}): M \mbox{ is an } m\times n\mbox{ matrix  with at most }
    s \mbox{ distinct entries}\} . $$
\end{abstract}

\section{Introduction}

The problem of studying the properties of matrices that avoid
certain submatrices or patterns is a classical and well studied
problem in combinatorics. It is investigated from a matrix point of
view as well as in an equivalent formulation of forbidden
subgraphs of bipartite graphs; see \cite{A},
\cite{F}, \cite{DRW}, \cite{PS}, et al.  Most of the previous research is
devoted to extremal and structural problems of matrices with no
forbidden submatrices. There are only a few results studying
efficient modifications of matrices or graphs such that the
resulting structure satisfies certain properties, for example,
\cite{EGR} and \cite{EGS}. In this paper, we apply powerful graph
theoretic techniques to study the distance properties between
certain classes of matrices. Our main goal is to investigate the
number of positions where the entry-changes need to be performed
on a given matrix such that the resulting matrix
does not have a fixed subpattern.
Although this problem is of independent theoretical interest, it
has multiple applications in computational biology such as in
the compatibility of evolutionary trees and in studying metabolic
networks, see \cite{CEFS}, \cite{SAN1998}.

For positive integers $m,n,s$, with  $s\leq mn$, let ${\cal
M}(m,n;s)$ denote the set of all $m\times n$ matrices with a fixed
number, $s$, of distinct entries.  Let $[m]\stackrel{\rm
def}{=}\{1,\ldots,m\}$.  Let ${\cal A}=\{ A_1,\ldots, A_r\}$ be a
partition of pairs from $[m]\times [n]$ into $r$ nonempty classes.
An $m\times n$ matrix $A=(a_{ij})$ is said to have a {\it pattern
${\cal A}$} provided that $a_{ij}=a_{i'j'}$ if and only if
$(i,j),(i',j')\in A_t$ for some $t\in\{1,\ldots,r\}$. It follows,
in particular, that two $m\times n$ matrices $A$ and $B$ with sets
of distinct entries $S(A)$ and $S(B)$, respectively, have the same
pattern if there is a bijection $g:S(A) \rightarrow S(B)$ such
that $B(i,j)= g(A(i,j))$ for all $1\leq i\leq m$ and all $1\leq
j\leq n$.

\begin{example}
Matrices $A$ and $B$ have the same pattern with a corresponding bijection $g$;
matrices $A$ and $B'$ have different patterns: \\
$$ A=\left(\begin{array}{ccc}
            1 & 4 & 3 \\
            1 & 1 & 4 \end{array}\right), \qquad
   B=\left(\begin{array}{ccc}
            5 & 1 & 2 \\
            5 & 5 & 1 \end{array}\right), \qquad
   B'=\left(\begin{array}{ccc}
            5 & 1 & 2 \\
            0 & 5 & 1 \end{array}\right). $$

\end{example}

In this case, $g(1)=5$, $g(4)=1$, $g(3)=2$. \\

A $k\times\ell$ matrix ${B}$ is a {\it submatrix} of an $m\times
n$ matrix  $A$ if there are nonempty subsets $\{i_1,\ldots,i_k\}$
and $\{j_1,\ldots,j_{\ell}\}$ of distinct indices with
$\{i_1,\ldots,i_k\}\subseteq [m]$,
$\{j_1,\ldots,j_{\ell}\}\subseteq [n]$ such that
$B(\alpha,\beta)=A(i_\alpha,j_{\beta})$, $1\leq \alpha \leq k$,
$1\leq \beta \leq l$.  If, for a matrix $M'$, there is a submatrix
$M$ with pattern ${\cal A}$ then we say that $M'$ has a {\it
subpattern ${\cal A}$}.

\begin{definition}
For a pattern ${\cal A}$  and positive integers $m,n,s$, we define
$\Forb(m,n; s, {\cal A})$ to be the set of all $m\times n$
matrices with at most $s$ distinct entries and not containing
subpattern ${\cal A}$.
\end{definition}

\begin{example} Let ${\cal
A}=\{\{(1,1),(1,2),(2,1)\},\{(2,2)\}\}$.  The set
$\Forb(m,n;2,{\cal A})$ consists of all $m\times n$ matrices which
have at most $2$ distinct entries and contain no submatrix of the
form $\left(\begin{array}{cc} x & x \\ x & y \end{array}\right)$,
$\left(\begin{array}{cc} y & x \\ x & x \end{array}\right)$,
$\left(\begin{array}{cc} x & x \\ y & x \end{array}\right)$,
$\left(\begin{array}{cc} x & y \\ x & x \end{array}\right)$,
$x\neq y$. In particular, $\Forb(m,n;2,{\cal A})$ consists of
$m\times n$ matrices with all entries equal  and all $m\times n$
matrices with two distinct entries such that each row has all
equal entries.
\end{example}

Next we define the distance between two matrices and between
classes of matrices.
For two matrices $A$ and $B$ of the same dimensions,  we say that
$\Dist(A,B)$ is the number of positions in which  $A$ and $B$ differ;
i.e., it is the matrix Hamming distance. For a class of matrices
${\cal F}$ and a matrix $A$, all of the same dimensions, we denote
$\Dist(A,{\cal F})= \min\{\Dist(A, F): F\in {\cal F}\}$. Finally,
$$f(m,n; s,{\cal A})=\max\{\Dist(A,{\cal F}):A\in{\cal M}(m,n; s), {\cal F}=\Forb(m,n; s, {\cal A}) \}.$$
This function
corresponds to the minimum number of positions on which the entries
need to be changed in any
$m\times n$ matrix with at most $s$ distinct entries in order to
eliminate all subpatterns ${\cal A}$.  This problem is also called
an {\it editing distance problem}, since we consider the minimum
number of editing operations on a matrix, where each editing
operation is a change of an entry in some position.

Note that $\Forb(m,n;s,{\cal A})$ might be  an empty set of
matrices for some patterns ${\cal A}$. For example, let $s$ be
fixed, and ${\cal A}$ be a pattern having exactly one set; i.e., a
pattern corresponding to matrices with all entries being equal.
We call such a  pattern a {\it trivial pattern}. If $m$ and $n$ are large,
then there is no $m\times n$ matrix with fixed number of distinct entries
avoiding pattern ${\cal A}$. This follows from the finiteness of the bipartite Ramsey number, see \cite{GRS}.
%This follows from the solution of the problem of Zarankiewicz
%\cite{Z} by F\"uredi, \cite{F}.
On the other hand, when a pattern ${\cal A}$ has at least two distinct entries, then the class
$\Forb(m,n;s,{\cal A})$ is nonempty since it contains all $m\times
n$ matrices with a trivial pattern. Our main result is the following:
\begin{theorem}
Let $s, r$ be positive integers, $s\geq r$.
Let $b_1,b_2$ be positive constants such that $b_1\leq m/n\leq
b_2$.  Let ${\cal A}$ be  a non-trivial  pattern with $r$ distinct
entries, then
$$ f(m,n; s,{\cal A})=
   (1+o(1))\left(\frac{s-r+1}{s}\right)mn . $$
   \label{matrix}
\end{theorem}

%Theorem  \ref{matrix} is an immediate consequence of
%Theorem~\ref{bipartite}.
%\begin{theorem} \renewcommand{\theenumi}{\alph{enumi}}
%Let $s, r$ be positive integers, $ s\geq r$.  Let $b_1, b_2$ be positive constants such that $b_1\leq m/n\leq
%b_2$.  Let ${\cal A}$ be a non-trivial  pattern with $r$ distinct entries.
%Then
%\begin{enumerate}
%   \item for any $m\times n$ matrix with at most  $s$ entries, it is
%   sufficient to change entries in at most
%   $(1+o(1))\left(\frac{s-r+1}{s}\right)mn$ positions such that
%   the resulting matrix does not contain a subpattern $\cal A$;
%   \item there is an $m\times n$ matrix $M$ with $s$ entries such
%   that in order to destroy all subpatterns ${\cal A}$ in $M$,
%   one needs to change the entries in at least
%   $(1+o(1))\left(\frac{s-r+1}{s}\right)mn$ positions.
%\end{enumerate} \label{bipartite}
%\end{theorem}

We shall prove these results using graph-theoretic formulations. A
graph $H=(V,E)$ is bipartite if its vertex set can be partitioned
such that $V=X\cup Y$, $X\cap Y=\emptyset$, and its edge set $E$
is a subset of $X\times Y$.  If $m=|X|$, $n=|Y|$ and $E=X\times
Y$, then this graph is denoted $K_{m,n}$ and called a complete
bipartite graph.  Now, we can introduce a pattern on the edges of
a complete bipartite graph as a partition of the edges in exactly the
same manner as above. Let ${\cal A} = \{A_1, \ldots, A_r\}$ such
that $ E = A_1 \cup \cdots \cup A_r$ and  $A_i$'s are  nonempty and
pairwise disjoint. Then ${\cal A}$ is called a pattern on $E$.
Now, let $c$ be a coloring of edges of $K_{m,n}$.
 We say that $c$ has a {\it pattern} ${\cal A}$ if it satisfies the property that
$c(e)=c(e')$ if and only if
$e,e'\in A_i$ for some $i=1,\ldots,r$. If $c$ is an edge-coloring
of a graph $G$, we say that a coloring $c'$ of a graph $G'$ occurs
in $G$ under coloring $c$ if there is a subgraph $H$ of $G$
isomorphic to $G'$ such that the coloring $c$ restricted to  $H$ coincides with
the coloring $c'$ of $G'$. Similar to the case with matrices,
for a color pattern ${\cal A}$ defined on the edges of a graph $G'$, we say that $G$ has a
subpattern ${\cal A}$ if there is an occurrence of a subgraph $H$
in $G$ such that $H$ is isomorphic to $G'$ and the coloring $c$ restricted to
$H$ has a pattern ${\cal A}$.

For two edge-colorings $c$ and $c'$ of a graph
$G$, we say that the {\it edit distance} between $c$ and $c'$ on
$G$ is the smallest number of edge-recolorings in $G$ colored
under $c$ needed to obtain $c'$.
For a given pattern ${\cal A}$ on edges of a complete bipartite
graph, and an edge-colored $K_{m,n}$ with coloring $c$, let
$F(m,n; c, {\cal A})$ be the smallest number of edge-recolorings of
$K_{m,n}$ colored by $c$ such that the resulting coloring does not contain a
subpattern ${\cal A}$. Define
$$F(m,n; s, {\cal A}) := \max \{F(m,n; c,
{\cal A}):  c \mbox { uses } s \mbox{ colors}\}.$$ \\

\noindent {\bf Observation.} There is a bijection $g$ between all
$m\times n$ matrices with $s$ distinct entries and all edge-colorings of
$K_{m,n}$ using $s$ colors. Indeed, this bijection can be defined
as $g(M(i,j))=c(\{i,j\})$, $i\in \{1, \ldots, m\}$, $j\in \{1,
\ldots, n\}$;   where $c(\{i,j\})$ is the color of an edge $\{i,j\}$
and  $M(i,j)$ is the $(i,j)$th entry of the  matrix.  Moreover, a
matrix $M$ does not have subpattern ${\cal A}$ if and only if a
coloring $g(M)$ does not have a subpattern ${\cal A}$. \\

For all other graph-theoretic terminology, we refer the reader to
\cite{W}. Our main theorem is proven in terms of graph colorings.

\begin{theorem}\label{graph}
Let  $\epsilon$,
$0<\epsilon<1$ be fixed, and let $m',n',s, r $;   $s\geq r$,  be fixed as well. Let $m+n$ be sufficiently large and let
$\cal A$ be a pattern on $K_{m',n'}$  with $r$ colors.
Then,
$$ \left(1-\epsilon\left(5s+2+(s+1)\left(\frac{m}{n}
                                         +\frac{n}{m}
                                   \right)\right)\right)
   \left(\frac{s-r+1}{s}\right)mn\leq F(m,n;s,{\cal A})
   \leq\left(\frac{s-r+1}{s}\right)mn . $$
\end{theorem}

Observe that now  Theorem \ref{matrix} is an immediate corollary of Theorem \ref{graph} which we
prove in Section 3. Section 2 describes the techniques that we use in the proof.

\section{Main tools}

For two disjoint sets of vertices $X$ and $Y$, we shall refer to a
pair $(X,Y)$  as a complete bipartite graph with partite sets $X$
and $Y$.  We denote its edges by
$E(X,Y)$. Let $c:E(X,Y) \rightarrow \{1, \ldots, s\}$ be an
edge-coloring of a pair $(X,Y)$. For each color $\nu\in \{1,
\ldots, s\}$, and any two subsets $X'\subseteq X$, $Y'\subseteq
Y$, we denote by $E_{\nu}(X',Y')$ the set of edges of color $\nu$
in a pair $(X',Y')$. Then $d_{\nu}(X', Y')$ is the {\it density of
a color $\nu$} in the subgraph induced by $X'$ and  $Y'$, defined
as follows:
$$ d_{\nu}(X',Y') = \frac{|E_{\nu}(X',Y')|}{|X'||Y'|}. $$
For $x\in X\cup Y$, we define $N_{\nu}(x)$ to be the set of all vertices joined to $x$ by edges of color $\nu$.
We say that a pair
$(X,Y)$ is {\it $\epsilon$-regular in color $\nu$} if for every
$X'\subseteq X$ and $Y'\subseteq Y$ with  sizes $|X'|\geq\epsilon
|X|$, $|Y'|\geq \epsilon |Y|$,  we have
\begin{equation}
\left|d_{\nu}(X,Y)-d_{\nu}(X',Y')\right|<\epsilon. \label{epreg}
\end{equation}

%\begin{lemma}[Many-Color Regularity Lemma~\cite{KS}]
%   For any $\epsilon>0$ and integers $r,m$ there exists
%   positive integers $\kappa$ and $M$ such that if the edges of a
%   graph $G$ are $r$-colored then the vertex set $V(G)$
%   can be partitioned into sets $V_0,V_1,\ldots,V_k$, for some
%   $m\leq k\leq M$, so that $|V_0|<\epsilon n$, $|V_i|=\kappa$
%   for every $i\geq 1$, and all but at most $\epsilon k^2$
%   pairs $(V_i,V_j)$ are $\epsilon$-regular in color $\nu$ for
%   each $\nu=1,\ldots,r$.
%   Moreover, if $V(G)$ is partitioned into sets $U_0, U_1, \ldots,
%   U_t$ such that all set $U_j$, $j=1, \ldots, t$ are of equal
%   size with $t\leq T$  and $|U_0|\leq (\epsilon/2) |V(G)|$ there
%   is a partition of $V(G)$ into sets $V_0,V_1,\ldots,V_k$
%   satisfying previously stated conditions with
%   $M=M(\epsilon,r,m,T)$ and the additional restriction that for
%   all $i=1, \ldots, k$ we have $V_i\subseteq U_j$ for some
%   $j\in \{1, \ldots, t\}$. \label{MCRL}
%\end{lemma}

%\begin{remark}

Lemma \ref{BMCRL} is based on the so-called  many-color  regularity lemma of
Szemer\'edi, see \cite{KSSSz},  and  is
an implication of the refinement argument, i.e., Theorem 8.4 in \cite{KS}.

%\end{remark}

\begin{lemma}[Bipartite Many-Color Regularity Lemma~\cite{KS}]
   For any $\epsilon>0$ and integers $s,m_0$ there exists
   $M$, a positive integer,  such that if the edges
   of a pair $(X,Y)$ are colored with  $1, \ldots, s$,
   then the vertex set
   $X\cup Y$ can be partitioned into sets $V_0,V_1,\ldots,V_k$,
   for some $k$,  $m_0\leq k\leq M$,  so
   that $|V_0|<\epsilon (|X|+|Y|)$, and $|V_i|=|V_j|$, for
   $i,j\in \{1,\ldots,k\}$, and all but at most $\epsilon k^2$
   pairs $(V_i,V_j)$ are $\epsilon$-regular in color $\nu$,
   for each $\nu=1,\ldots,s$,
   and either $V_i\subseteq X$ or $V_i\subseteq Y$, for
   $i=1,\ldots,k$. \label{BMCRL}
\end{lemma}

%\begin{proof}
%We apply lemma \ref{BMCRL} with a partition $U_0, \ldots, U_t$
%such that for each $i=1, \ldots, t$, either $U_i\subseteq X$ or
%$U_i\subseteq Y$ and $|U_0|\leq (\epsilon/2)(|X|+|Y|)$.
%\end{proof}

In addition, we need to prove a multicolor version of the
so-called intersection property, which is stated in~\cite{KS} and revised
in~\cite{AKM}.

\begin{fact}[Many-Color Intersection Property]
Let $\epsilon>0$ and $\delta>0$  be fixed and $r$ and $\ell$ be
positive integers. Let $(A,B)$ be a pair with edges colored such
that color $\nu$ is $\epsilon$-regular with density $d_{\nu}$,
$ d_{\nu}\geq \delta$ for $\nu=1,\ldots,r$. Let
$Y\subset B$. Assume that $(\delta-\epsilon)^{\ell-1}|Y|>\epsilon
|B|$. Let $k_{\nu}$ for $\nu=1,\ldots,r$ be a positive integer
such that $\sum_{\nu=1}^rk_{\nu}=\ell$ and let any vector ${\bf
a}\in A^{\ell}$ be indexed such that
$$ {\bf a}=\left(a_{[1,1]},\ldots,a_{[1,k_1]},
   a_{[2,1]},\ldots,a_{[r-1,k_{r-1}]},
   a_{[r,1]},\ldots,a_{[r,k_{r}]}\right) . $$
Then,
\begin{equation}
   \#\left\{{\bf a}\in A^{\ell} :
   \left|Y\cap\bigcap_{\nu=1}^{r}\bigcap_{i=1}^{k_{\nu}}
   N_{\nu}(a_{[\nu,i]})\right|
   <\prod_{\nu=1}^r\left(d_{\nu}-\epsilon\right)^{k_{\nu}}|Y|\right\} \leq \ell\epsilon |A|^{\ell}. \label{lb}
  % \#\left\{{\bf a}\in A^{\ell} :
  % \left|Y\cap\bigcap_{\nu=1}^{s}\bigcap_{i=1}^{k_{\nu}}
  % N_{\nu}(a_{[\nu,i]})\right|
  % >\prod_{\nu=1}^s\left(d_{\nu}+\epsilon\right)^{k_{\nu}}|Y|\right\}
  % & \leq \ell\epsilon |A|^{\ell} . & \label{ub}
\end{equation}
\label{fRIP}
\end{fact}

The proof of Fact \ref{fRIP} is a standard argument which  follows by induction on $\ell$.

\begin{corollary}\label{any}
Let $\epsilon>0$ and $\delta>0$  be  fixed  and  $r$  and  $\ell$ be positive integers.
 Let $c$ be an edge-coloring of a pair  $(A,B)$ with  at least $r$ colors from  $\{1, \ldots, r, \ldots \}$
such that color $\nu$ is
$\epsilon$-regular with density $d_{\nu}$, $ d_{\nu}\geq \delta$, for
$\nu=1,\ldots, r$. Let us be
given that $(\delta-\epsilon)^{\ell-1}>\epsilon$ and
$2r^{\ell}\ell\epsilon<1$ and $(\delta-\epsilon)^{\ell}|B|\geq \ell$.
Then any  edge-coloring of $K_{\ell,\ell}$ with colors from
$\{1,\ldots,r\}$ will occur as a subcoloring of $c$.
%In particular, if we have $\ell\geq 2$, $\delta\geq 2\epsilon$,
%and $\epsilon<\ell^{-1}s^{-\ell}/2$, then each coloring of
%$K_{\ell,\ell}$ occurs in $(A,B)$. \label{cRIP}
\end{corollary}

\section{Proof of Theorem \ref{graph}}
\subsection{Upper bound}
We shall show that for any $s$-edge-coloring of a complete
bipartite graph with vertex class of sizes $m$ and $n$, there are
at most $\left(\frac{s-r+1}{s}\right)mn$ editing operations
sufficient  to destroy a fixed color pattern with $r$ colors.

Let $\cal A$ be a color pattern with $r$ sets defined on a complete
bipartite graph $G$  and let $c$ be an edge-coloring of $K_{m,n}$ with
$s$ colors. Without loss of generality, let $1$ be the color of
the largest color class in $c$. We shall recolor the $s-r+1$ smallest color
classes of $c$ so that their new color is $1$. The resulting coloring will use only $r-1$ colors
and thus will not contain a forbidden pattern. The $s-r+1$
smallest color classes account for at most
$\left(1-(r-1)/s\right)mn$ edges.  Thus,
$$ F(n,m;s, {\cal A} )\leq\left(\frac{s-r+1}{s}\right)mn . $$

\subsection{Lower bound}
To establish the lower bound, we show that there is a coloring of the given complete
bipartite graph requiring many edit-operations to destroy a forbidden pattern.
We begin with a claim that gives us a coloring which is highly
regular.

\begin{claim}
Let $s$ be a positive integer, and  $0<\epsilon<1/2$. There
is an  integer $M$ such that if $|X|\geq M$ and $|Y|\geq
M$ then  there is an edge-coloring $c$ of a complete bipartite
graph $G=X\times Y$, with colors $1, 2, \ldots, s$, satisfying the following
property: If $X'\subseteq X$ and $Y' \subseteq Y$, such that
$|X'|,|Y'|>(|X|+|Y|)(1-\epsilon)/M$, then $d_{\nu}(X',Y')\in
(1/s-\epsilon,1/s+\epsilon)$, $\nu = 1, \ldots, s$. \label{rand}
\end{claim}

Claim 1 follows from standard applications of the Chernoff bound
(see~\cite {JLR}, Chapter 2).

%\noindent {\it Proof of Claim.} Consider a coloring formed by
%coloring each edge independently at random with $k$ colors so that
%the probability that an edge is colored with a color $i$ is $1/k$,
%for $i=1, \ldots, k$.
%\begin{eqnarray*}
%   \hspace{-1in}\lefteqn{\Pr\left\{\exists X',Y',\nu :
%   \quad X'\subseteq X, Y'\subseteq Y,\quad
%   |X'|,|Y'|>(|X|+|Y|)/2M \mbox{ and }
%   \left|d_{\nu}(X',Y')-1/k\right|\geq\epsilon\right\}}
%   \hspace{1in} \\
%   & \leq &
%   2 \left(2^{|X|+|Y|}\right)^2k\Pr\left\{d_{\nu}(X',Y')\leq
%   1/k-\epsilon\right\} \\
%   & \leq & 2k\cdot 2^{2(|X|+|Y|)}\cdot
%   \exp\left(-2(\epsilon|X'||Y'|)^2/(|X'||Y'|)\right) \\
%   & \leq & 2k\cdot 2^{2(|X|+|Y|)}\cdot
%   \exp\left(-2\epsilon^2|X'||Y'|\right) \\
%   & \leq & 2k\cdot 2^{2(|X|+|Y|)}\cdot
%   \exp\left(-2\epsilon^2\left((|X|+|Y|)/(2M)\right)^2\right)
%   <1 ,
%\end{eqnarray*}
%if $|X|+|Y|$ is large enough. Here, we use a Chernoff bound, see
%for example \cite{JLR}. Therefore, there is at least one coloring,
%$c$, having the property that all pairs of sufficiently large sets
%have density within $\epsilon$ of the expectation $1/k$ in each
%color. This proves the claim.\\ \\

Fix $\epsilon>0$, let $c'$ be a coloring of the pair $(X,Y)$
$|X|=m, |Y|=n$,  of minimum edit distance from $c$ with the
property that $c'$ contains no subpattern $\cal A$. Apply
Lemma~\ref{BMCRL} with parameters $\epsilon$, $s$ and $m_0=1$ to
the coloring $c'$. Let $M$ be the constant given by
Lemma~\ref{BMCRL} and the partition having all the non-leftover
sets being enumerated as $X_1,\ldots,X_{p},Y_1,\ldots,Y_{q}$ with
$|X_i|=|Y_j|=Q$  and $X_i\subseteq X$, $Y_j \subseteq Y$ for
$1\leq i\leq p$, $1\leq j\leq q$. We call a pair $(X_i, Y_j)$, a
{\it good} pair if it is $\epsilon$-regular in each color $\nu \in
\{1, 2, \ldots, s\}$ in coloring $c'$. We have that there are at
most  $s\epsilon(p+q)^2$ pairs which are not good. Moreover, for
each good pair $(X_i, Y_j)$ there are at most $r-1$ colors such
that the denstity of those classes in coloring $c'$ is at least
$\delta=2\epsilon$. Otherwise, Corollary~\ref{any} would imply
that pattern ${\cal A}$ appears in $c'$, a contradiction.
Therefore, for a good  pair $(X_i,Y_j)$, there are at least
$(s-r+1)\left(\frac{1}{s}-3\epsilon\right)Q^2$ editing operations
needed to obtain coloring $c'$ from the coloring $c$. The
regularity lemma gives that $m\geq pQ\geq m-\epsilon(m+n)$ and
$n\geq qQ\geq n-\epsilon(m+n)$.  Therefore, the total number of
recolored edges is at least
\begin{eqnarray*}
   \lefteqn{(s-r+1)\left(\frac{1}{s}-3\epsilon\right)Q^2
            \left(pq-s\epsilon(p+q)^2\right)} \\
   & \geq & \left(\frac{s-r+1}{s}\right)(1-3s\epsilon)
            \left(pQqQ-s\epsilon(pQ+qQ)^2\right) \\
   & \geq & \left(\frac{s-r+1}{s}\right)(1-3s\epsilon)
   \left(\left(m-\epsilon(m+n)\right)\left(n-\epsilon(m+n)\right)
   -s\epsilon(m+n)^2\right) \\
%   & \geq & \left(\frac{s-r+1}{s}\right)mn(1-3s\epsilon)
%   \left(\left(1-\epsilon\left(1+\frac{n}{m}\right)\right)
%         \left(1-\epsilon\left(\frac{m}{n}+1\right)\right)
%         -s\epsilon\left(1+\frac{n}{m}\right)
%                   \left(\frac{m}{n}+1\right)\right) \\
%   & \geq & \left(\frac{s-r+1}{s}\right)mn(1-3s\epsilon)
%   \left(1-\epsilon\left(2+\frac{m}{n}+\frac{n}{m}\right)
%          -s\epsilon\left(1+\frac{m}{n}\right)
%                    \left(1+\frac{n}{m}\right)\right) \\
%   & \geq & \left(\frac{s-r+1}{s}\right)mn
%            \left(1-\epsilon
%                    \left(3s+2+\frac{m}{n}+\frac{n}{m}
%                          +s\left(1+\frac{m}{n}\right)
%                            \left(1+\frac{n}{m}\right)
%                    \right)\right) \\
  & \geq & \left(\frac{s-r+1}{s}\right)mn
            \left(1-\epsilon
                    \left(5s+2
                          +(s+1)\left(\frac{m}{n}
                                      +\frac{n}{m}
                                \right)\right)\right) .
%             \\
%    & \geq & \left(\frac{s-r+1}{s}\right)mn(1-C\epsilon) ,
\end{eqnarray*}
% where $C=5s+2+(s+1)(b_1+b_2)$.
\proofend

\noindent
{\bf Remark.} It should be noted that, although we prove theorems for
  submatrices, our results easily follow for other patterns.
  Suppose we wish to forbid patterns of the form
  $\left(\begin{array}{rr} 1 & 2 \\ 1 & * \end{array}\right)$,
  where the  $*$ represents any entry, either a repeated $1$ or
  $2$ or a new entry $3$.  Our result depends only on the number
  of distinct entries in the pattern, so the (asymptotic) number
  of changes necessary and sufficient to forbid this pattern is
  the same as the number of changes needed to forbid $\left(\begin{array}{rr} 1 & 2 \\ 1 & 1
  \end{array}\right)$ or $\left(\begin{array}{rr} 1 & 2 \\ 1 & 2
  \end{array}\right)$ (that is,
  $(1+o(1))\left(\frac{s-1}{s}\right)mn$) but fewer than to forbid
  $\left(\begin{array}{rr} 1 & 2 \\ 1 & 3 \end{array}\right)$
  (that is,  $(1+o(1))\left(\frac{s-2}{s}\right)mn$).\\

\noindent
{\bf Acknowledgement.} We are indebted to anonymous referees whose careful reading
and friendly suggestions helped to significantly improve the presentation of the results.

\end{document}